\documentclass[a4paper,12pt]{article}
\begin{document}

\title{Analytic extensions of functions of Cayley-Dickson variables.}
\author{Sergey V. Ludkowski}
\date{23 April 2018}
\maketitle
\begin{abstract}
Analytic approximations of functions of Cayley-Dickson variables are
investigated. The case of functions of complexified Cayley-Dickson
variables is also encompassed. Moreover, extensions of functions of
Cayley-Dickson variables are studied. \footnote{key words and
phrases: function; analytic; extension; 
Cayley-Dickson algebra  }

\end{abstract}
\par Address: Dep. Appl. Mathematics, Moscow State Techn. Univ. MIREA,
\par av. Vernadsky 78, Moscow 119454, Russia; e-mail: sludkowski@mail.ru

\section{Introduction.}
The theory of approximations and extensions of differentiable
functions of real or complex variables is important for studies of
real or complex manifolds and boundary value problems of PDEs (see,
for example,
\cite{henlei}-\cite{hormbl,kling,kodai,michor,svalkorplb,zabokoportj16}
and references therein). \par On the other hand, since the years
2000-th hypercomplex analysis over Cayley-Dickson algebras began to
be developed. In particular, it includes analysis over octonions. It
has found many applications in mathematical analysis, partial
differential equations, mathematical physics, noncommutative
geometry and natural sciences (see \cite{baez,boloktodb},
\cite{dickson,frenludkfejms18}-\cite{guetze,kansol,ludkrimut14,ludfov}-\cite{ludkcvee13}
and references therein). Founders of quantum field theory such as
Yang and Mills formulated a problem of developing analysis over
octonions and Cayley-Dickson algebras for studying physical
problems.  As it was shown in
\cite{ludkaaca2014,ludfov}-\cite{ludkfscdvjms07} a differentiability
of functions in terms of the word algebra \cite{bourbalg,razm}
provides abundant families of functions of octonion and
Cayley-Dickson variables in comparison with the classical
differentiability.
\par Extensions of differentiable functions of real variables
were described in the paper \cite{whitntrams34}.
\par In this article analytic approximations of functions of Cayley-Dickson variables are
investigated. The case of functions of complexified Cayley-Dickson
variables is also encompassed. Moreover, extensions of functions of
Cayley-Dickson variables are studied. The main results are Theorems
6-8 about analytic extensions.
\par The main results of this paper are obtained for the first time.

\section{Analytic approximations and extensions.}
For avoiding misunderstandings we first remind definitions and
notations. A reader familiar with octonions and Cayley-Dickson
algebras may skip Notes 0 and 1.
\par {\bf 0. Note.} The algebra $\bf O$ of octonions (octaves, the Cayley algebra)
is defined as an eight-dimensional nonassociative algebra over $\bf
R$ with a basis, for example,
\par $(1)$ ${\bf b}_3:={\bf b}:= \{ 1, i, j, k, l, il, jl, kl \} $ such that
\par $(2)$ $i^2=j^2=k^2=l^2=-1$, $ij=k$, $ji=-k$,
$jk=i$, $kj=-i$, $ki=j$, $ik=-j$, $li=-il$, $jl=-lj$, $kl=-lk$;
\par $(3)$ $(\alpha +\beta l)(\gamma +\delta l)=(\alpha \gamma
-{\tilde {\delta }}\beta )+(\delta \alpha +\beta {\tilde {\gamma }})l$ \\
is the multiplication law in the octonion algebra $\bf  O$ for each
$\alpha $, $\beta $, $\gamma $, $\delta \in \bf H$, where $\bf H$
denotes the quaternion skew field, $\xi :=\alpha +\beta l\in \bf O$,
$\eta :=\gamma +\delta l\in \bf O$, ${\tilde z}:=v-wi-xj-yk$ for
$z=v+wi+xj+yk \in \bf H$ with $v, w, x, y\in \bf R$.
\par The octonion algebra is neither commutative, nor associative,
since $(ij)l=kl$, $i(jl)=-kl$, but it is distributive and ${\bf R}1$
is its center. If $\xi :=\alpha +\beta l\in \bf O$, then
\par $(4)$ ${\tilde {\xi }}:={\tilde {\alpha }}-\beta l$ is called
the adjoint element of $\xi $, where $\alpha , \beta \in \bf H$.
Then
\par $(5)$ $(\xi \eta )^{\tilde .}={\tilde {\eta }}{\tilde {\xi }}$,
${\tilde {\xi }} + {\tilde {\eta }}= (\xi +\eta )^{\tilde .}$ and
$\xi {\tilde {\xi }}=|\alpha |^2+|\beta |^2$, \\
where $|\alpha |^2=\alpha {\tilde {\alpha }}$ such that
\par $(6)$ $\xi {\tilde {\xi }}=:|\xi |^2$ and $|\xi |$
is the norm in $\bf O$. Therefore,
\par $(7)$ $|\xi \eta |=|\xi | |\eta |$, \\
consequently, $\bf O$ does not contain divisors of zero (see also
\cite{kansol,kurosh,waerd,ward}). The multiplication of octonions
satisfies equations $(8,9)$ below:
\par $(8)$ $(\xi \eta )\eta =\xi (\eta \eta )$,
\par $(9)$ $\xi (\xi \eta )=(\xi \xi )\eta $, \\
that forms the alternative system. In particular, $(\xi \xi )\xi
=\xi (\xi \xi )$. Put ${\tilde {\xi }}=2a-\xi $, where $a=Re (\xi
):=(\xi + {\tilde {\xi }})/2 \in \bf R$. Since ${\bf R}1$ is the
center of the octonion algebra $\bf O$ and ${\tilde {\xi }}\xi =\xi
{\tilde {\xi }}=|\xi |^2$, then from $(8,9)$ by induction it
follows, that for each $\xi \in \bf O$ and each $n$-tuplet
(product), $n\in \bf N$, $\xi (\xi (... \xi \xi )... )=(...(\xi \xi
)\xi ...)\xi $ the result does not depend on an order of brackets
(order of multiplications), hence the definition of $\xi ^n:= \xi
(\xi (... \xi \xi )...)$ does not depend on the order of brackets.
This also shows that $\xi ^m\xi ^n=\xi ^n\xi ^m$, $\xi ^m{\tilde
{\xi }}^m={\tilde {\xi }}^m\xi ^n$ for each $n, m\in \bf N$ and $\xi
\in \bf O$.
\par Apart from the quaternions, the octonion algebra can not be
realized as the subalgebra of the algebra ${\bf M}_8({\bf R})$ of
all $8\times 8$-matrices over $\bf R$, since $\bf O$ is not
associative, but the matrix algebra ${\bf M}_8({\bf R})$ is
associative (see, for example, \cite{baez,kurosh,waerd,ward}). There
are the natural embeddings ${\bf C}\hookrightarrow \bf O$ and ${\bf
H}\hookrightarrow \bf O$, but neither $\bf O$ over $\bf C$, nor $\bf
O$ over $\bf H$, nor $\bf H$ over $\bf C$ are algebras, since the
centers of them are $Z({\bf H})=Z({\bf O})=\bf R$ equal to the real
field.
\par We consider also the Cayley-Dickson algebras ${\cal A}_n$
over $\bf R$, where $2^n$ is its dimension over $\bf R$. They are
constructed by induction starting from $\bf R$ such that ${\cal
A}_{n+1}$ is obtained from ${\cal A}_n$ with the help of the
doubling procedure, in particular, ${\cal A}_0:=\bf R$, ${\cal
A}_1=\bf C$, ${\cal A}_2=\bf H$, ${\cal A}_3=\bf O$ and ${\cal A}_4$
is known as the sedenion algebra \cite{baez}. The Cayley-Dickson
algebras are $*$-algebras, that is, there is a real-linear mapping
${\cal A}_n\ni a\mapsto a^*\in {\cal A}_n$ such that
\par $(10)$ $a^{**}=a$,
\par $(11)$ $(ab)^*=b^*a^*$ for each $a, b\in {\cal A}_n$.
Then they are nicely normed, that is,
\par $(12)$ $a+a^*=:2 Re (a) \in \bf R$ and
\par $(13)$ $aa^*=a^*a>0$ for each $0\ne a\in {\cal A}_n$.
The norm in it is defined by the equation:
\par $(14)$ $|a|^2:=aa^*$.
\par We also denote $a^*$ by $\tilde a$. Each non-zero Cayley-Dickson
number $0\ne a\in {\cal A}_n$ has a multiplicative inverse given by
$a^{-1}=a^*/|a|^2$.
\par The doubling procedure is as follows. Each $z\in {\cal A}_{n+1}$
is written in the form $z=a+bl$, where $l^2=-1$, $l\notin {\cal
A}_n$, $a, b \in {\cal A}_n$. The addition is componentwise. The
conjugate of a Cayley-Dickson number $z$ is prescribed by the
formula:
\par $(15)$ $z^*:=a^* -bl$. \\
The multiplication is given by:
\par $(16)$ $(a +bl)(c +dl)=(a c -{\tilde {d}}b)+(d a + b {\tilde {c}})l$\\
for each $a, b, c, d$ in ${\cal A}_n$.

\par {\bf 1. Note.} By $ \{ i_0, i_1,...,i_{2^r-1} \} $  is denoted the standard basis
of the Cayley-Dickson algebra ${\cal A}_r={\cal A}_{r,\bf R}$ over
the real field $\bf R$ such that $i_0=1$, $i_l^2=-1$ and
$i_li_k=-i_ki_l$ for each $l\ne k$ with $1\le l$ and $1\le k$. For
$r\ge 3$ the multiplication of them is generally nonassociative (see
also Note 0). In particular ${\cal A}_3$ is the octonion algebra.
Henceforth the complexified Cayley-Dickson algebra ${\cal A}_{r,\bf
C} = {\cal A}_{r} \oplus ({\cal A}_{r}{\bf i})$ also is considered,
where ${\bf i}^2=-1$, $~ {\bf i}b=b{\bf i}$ for each $b\in {\cal
A}_r$, $~2\le r<\infty $. This means that each complexified
Cayley-Dickson number $z\in {\cal A}_{r,C}$ can be written in the
form $z=x+{\bf i}y$ with $x$ and $y$ in ${\cal A}_r$,
$x=x_0i_0+x_1i_1+...+x_{2^r-1}i_{2^r-1}$, while $x_0,...,x_{2^r-1}$
are in $\bf R$. The real part of $z$ is $Re (z)=x_0=(z+z^*)/2$, the
imaginary part of $z$ is defined as $Im (z)=z-Re (z)$, where the
conjugate of $z$ is $z^*=\tilde{z} = Re (z) - Im (z)$. Thus
$z^*=x^*-{\bf i}y$ with $x^*=x_0i_0-x_1i_1-...-x_{2^r-1}i_{2^r-1}$.
Then $|z|^2=|x|^2+|y|^2$, where $|x|^2=xx^*=x_0^2+...+x_{2^r-1}^2$.
\par If $U$ is a domain in ${\bf F}^{l2^r}$, then to
to each vector $u=(u_0,...,u_{2^r-1}) \in U$ a unique Cayley-Dickson
number $z=\hat{z}(u)=u_0i_0+u_1i_1+...+u_{2^r-1}i_{2^r-1}$ is posed,
where either ${\bf F}={\bf R}$ or ${\bf F}={\bf C}$, $~u_j\in {\bf
F}^l$ for each $j=0,...,2^r-1$; $~l\in {\bf N}$. This gives a domain
$V = \{ z: ~ z=\hat{z}(u), ~ u\in U \} $ in ${\cal A}_{r,\bf F}^l$.
Vise versa to each Cayley-Dickson number $z\in V$ a unique vector
$\pi (z)=\hat{u}(z)=u\in U$ corresponds:
\par $(1.1)$ $u_j=\pi _{j}(z)$ for each $j$, where $\pi _j : {\cal
A}_{r,\bf F}^l\to {\bf F}^l$ is a $\bf F$-linear operator given by
the formulas:
\par $(1.2)$ $\pi _j(z)=(-zi_j+ i_j(2^r-2)^{-1} \{ -z
+\sum_{k=1}^{2^r-1}i_k(zi_k^*) \} )/2$ for each $j=1,2,...,2^r-1$,
\par $(1.3)$ $\pi _0(z)=(z+ (2^r-2)^{-1} \{ -z +
\sum_{k=1}^{2^r-1}i_k(zi_k^*) \} )/2$, \\
where $ 2\le r\in \bf N$, $z$ is a Cayley-Dickson vector (or a
number for $l=1$) presented as
\par $(1.4)$ $z=z_0i_0+z_1i_1+...+z_{2^r-1}i_{2^r-1}\in {\cal
A}_{r,\bf F}$, $z_j\in {\bf F}^l$ for each $j=0,...,2^r-1$;
$z^*=z_0i_0-z_1i_1-...-z_{2^r-1}i_{2^r-1}$ (see Formulas
II$(1.1)-(1.3)$ in \cite{ludifeqcdla}). \par In the case ${\bf
F}={\bf R}$ by an ${\cal A}_r$-analytic function $f: V\to {\cal
A}_r$ on a domain $V$ open in ${\cal A}_r^l$ is meant a locally
analytic function $f$ in a $z$-representation. That is for each
marked point $\xi \in V$ an open neighborhood $V_{\xi }$ exists with
a phrase $\mu =\mu _{\xi }$ on $V_{\xi }$ such that $f=ev (\mu )$ on
$V_{\xi }$ and $D_{z^*}f(z)=0$ on $V$, where $ev$ denotes an
evaluation map of the phrase $\mu $ of a word algebra over ${\cal
A}_r^l$ (see Subsection 2.4 in \cite{ludkaaca2014} or
\cite{ludfov}).
\par If ${\bf F}={\bf C}$ we consider $f(z): W\to {\cal A}_{r,C}$ as a
function $g(x,y)=f(z)$ with $z=x+{\bf i}y$, where $z\in W\subset
{\cal A}_{r,C}^l$, $x$ and $y$ belong to ${\cal A}_r^l$, $g=\mbox{
}^fg$, $g=\mbox{ }_0g+{\bf i} \mbox{ }_1g$ with $\mbox{ }_vg\in
{\cal A}_r$ for $v=0$ and $v=1$. To an open set $W$ in ${\cal
A}_{r,C}^l$ an open set $V$ in ${\cal A}_r^{2l}$ is posed. If
$g(x,y)$ is analytic in $(x,y)$ on $V$, that is $\mbox{ }_vg\in
{\cal A}_r$ is ${\cal A}_r$-analytic on $V$ for $v=0$ and $v=1$,
then $f(z)$ will be called ${\cal A}_r$-analytic on $W$. \par
Therefore, functions $f$ on domains $V$ in ${\cal A}_r^l$ will be
below considered with values in ${\cal A}_{r,\bf F}$ if something
another will not be specified.
\par By a $C^m$ function $f: V\to {\cal A}_{r,\bf F}$ on an open domain or a canonically
closed domain $V$ in ${\cal A}_r^l$ is meant a $m$ times
continuously differentiable function in the corresponding real
variables $u$ on $U=\pi (V)$, where $0\le m \le \infty $. A family
of all $C^m$ functions on $V$ with values in ${\cal A}_{r,\bf F}$ is
denoted by $C^m(V,{\cal A}_{r,\bf F})$. For a function $G\in
C^m(V,{\cal A}_{r,\bf F})$ and each $|k|\le m$ the shortened
notation $\partial ^k_zG(z)$ of partial derivatives $\frac{\partial
^{|k|} G(z)} {\partial \mbox{}_1z_0^{k(1,0)} ...\partial
\mbox{}_lz_0^{k(l,2^r-1)}}$ will be used, where
\par $k=(k(1,0),...,k(1,2^r-1),...,k(l,0),...,k(l,2^r-1))$;
\par  $|k|=k(1,0)+...+k(l,2^r-1)$; $ ~ k(p,j)\in {\bf N}_0= \{ 0, 1, 2,... \} $
for each $p$ and $j$; where $\partial ^0_zG(z)=G(z)$ for a
unification of the notation.
\par The prefix "super" for superdifferentiable functions $f$ of ${\cal
A}_r$ variables will be omitted for brevity as in
\cite{ludkaaca2014}. If $D$ is a differentiation operator, by
$\frac{df(z)}{dz}$ or $D_zf(z)$ is denoted a derivative operator in
the $z$-variable on $V$ of an ${\cal A}_r$-differentiable function
$f: V\to  {\cal A}_{r,\bf F}$; $ ~ \frac{df(z)}{dz}.v$ or
$(D_zf(z)).v$ denote a derivative of $f$ along $v\in {\cal A}_r^l$.

\par {\bf 2. Definition.} Let $V$ be a set contained in ${\cal A}_r^l$.
Suppose that on $V$ functions $f_k(z): V\to {\cal A}_{r,\bf F}$ are
defined for each $|k|\le m$, where
\par $k=(k(1,0),...,k(1,2^r-1),...,k(l,0),...,k(l,2^r-1))$, \\ $0\le
k_{p,j}\in {\bf Z}$ for each $p$ and $j$, where
\par $|k|=k(1,0)+...+k(1,2^r-1)+...+k(l,0)+...+k(l,2^r-1)$. \par Let also
these functions satisfy the conditions: $$(2.1)\quad f_k(y)=
\sum_{s; ~ |s|\le m-|k|} \frac{f_{k+s}(z)}{s!}(\pi (y-z))^s +
Y_k(y;z)$$ for each $|k|\le m$ and each $y$ and $z$ in $V$, where
$y=(\mbox{}_1y,...,\mbox{}_ly)\in {\cal A}_r^l$, $\mbox{}_py\in
{\cal A}_r$ for each $p=1,...,l$; $\pi (y)=(\pi (\mbox{}_1y),..,\pi
(\mbox{}_ly))$; $\pi (\mbox{}_py)=(\pi _0(\mbox{}_py),...,\pi
_{2^r-1}(\mbox{}_py));$
$$u^k=\prod_{p,j}u_{p,j}^{k(p,j)};\quad s!= \prod_{p,j}s_{p,j}!$$
$u=(u_{1,0},...,u_{1,2^r-1},...,u_{l,0},...,u_{l,2^r-1})\in {\bf
R}^{l2^r}$; $ ~ (u_{p,0},...,u_{p,2^r-1})=\pi (\mbox{}_py)$ (see
also Formulas $(1.1)$-$(1.3)$). Functions $Y_k$ are supposed to be
satisfying the conditions: for each $z\in V$ and $\epsilon >0$ a
positive number $\delta >0$ exists such that for every $x$ and $y$
in $V$ with $|x-z|<\delta $ and $|y-z|<\delta $ and $|k|\le m$
\par $(2.2) \quad |Y_k(x;y)|\le |x-y|^{m-|k|}\epsilon $.
\par If conditions $(2.1)$ and $(2.2)$ are fulfilled, then it will be
said that the function $f(z)=f_0(z)$ is of class $C^m$ on $V$ in
terms of the functions \par $\{ f_k(z): ~ k\in {\bf N}_0^n, ~ |k|\le
m \} $, where ${\bf N}_0= \{ 0,  1, 2, ... \} $, $n=l2^r$.
\par If this is satisfied for each $m\in {\bf N} = \{ 1, 2, 3, ... \} $, then $f$ is of
class $C^{\infty }$ on $V$ in terms of the functions $\{ f_k(z): ~
k\in {\bf N}_0^n \} $.
\par By $C^{-1}(V,{\cal A}_{r,\bf F})$ is denoted a family of all
functions from $V$ into ${\cal A}_{r,\bf F}$.
\par {\bf 3. Lemma.} {\it Let $g: V\to {\cal A}_{r,\bf F}$ be of class
$C^m$ on a bounded closed subset $V$ in ${\cal A}_r^l$ in terms of
functions $\{ g_k(z):  ~ k\in {\bf N}_0^{l2^r}, ~  |k|\le m \} $,
where $m\in {\bf N}_0$, $ ~ 3\le r<\infty $, $~l\in {\bf N}$. Then
for each $\epsilon >0$ an ${\cal A}_r$-analytic function $G$ on
${\cal A}_r^l$ exists so that for each $|k|\le m$
$$(3.1)\quad |\partial ^k_z[G(z)- g(z)]|<\epsilon \mbox{ on } V.$$}
\par {\bf Proof.} Using formulas $(1.2)$ and $(1.3)$ for each
$z=(\mbox{}_1z,...,\mbox{}_lz)\in {\cal A}_r^l$ we write
$|\mbox{}_pz|^2$ and $|z|^2=\sum_{p=1}^l|\mbox{}_pz|^2$ in the
$z$-representation, where $\mbox{ }_pz\in {\cal A}_r$ for each $p$.
Hence to a function $f(z)=\exp (\sum_{p=1}^l b_p |\mbox{ }_pz|^2)$ a
phrase in the $z$-representation is posed
$$(3.2)\quad \mu = \exp (-(2^r-2)^{-1} \sum_{p=1}^l b_p [\mbox{}_pz
\sum_{k=0}^{2^r-1}i_k(\mbox{}_pzi_k)]),$$  where $b_p$ is a real
constant for each $p$. Thus $f(z-w)=ev_w(\mu )(z)$ for each $z\in
{\cal A}_r^l$ and $w\in {\cal A}_r^l$, where $ev_w$ means a
valuation map for a marked $w$. In particular for every marked
Cayley-Dickson number $w\in {\cal A}_r^l$ and a negative constant
$b_p$ for each $p=1,...,l$ the function $f(z-w)=ev_w(\mu )(z)$ is
positive, ${\cal A}_r$-differentiable bounded on ${\cal A}_r^l$ and
tends to zero when $|z|$ tends to the infinity.
\par To a differential form $d\mbox{}_1z_0\wedge ... \wedge d\mbox{}_1z_{2^r-1}
\wedge ... \wedge d\mbox{}_lz_0\wedge ... \wedge
d\mbox{}_lz_{2^r-1}$ on the Euclidean space ${\bf R}^{l2^r}$ a
differential form \par $(3.3)$ $d\pi _0(\mbox{}_1z)\wedge ... \wedge
d\pi _{2^r-1}(\mbox{}_1z) \wedge ... \wedge d\pi
_0(\mbox{}_lz)\wedge ... \wedge d\pi _{2^r-1}(\mbox{}_lz)=:\omega
(dz)$ \\ on ${\cal A}_r^l$ corresponds, where $d\pi _j(\mbox{}_kz)=
\pi _j(d \mbox{ }_kz)$ and $d \mbox{ }_kz=(d \mbox{
}_kz_0)i_0+...+(d \mbox{ }_kz_{2^r-1})i_{2^r-1}$ for each $j$ and
$k$. By $B({\cal A}_r^l,w,\rho )$ a closed ball in ${\cal A}_r^l$
with center at $w$ and of radius $\rho >0$ will be denoted,
\par $B({\cal A}_r^l,w,\rho )= \{ z\in {\cal A}_r^l: ~ |z-w|\le \rho
\} $. Then we put
$$(3.4) \quad \Phi (\rho ) = T \int_{B({\cal A}_r^l,0,\rho )}
\exp ( \sum_{p=1}^l  [\mbox{}_pz
\sum_{k=0}^{2^r-1}i_k(\mbox{}_pzi_k)])\omega (dz)$$ for each $0\le
\rho <\infty $ and choose a positive constant $T>0$ so that
\par $\lim_{\rho \to \infty } \Phi (\rho )=1$.
\par Let $0<\delta _0<\infty $ and let $V_1$ be a neighborhood of $V$ such that
$V_1\subset V^{\delta _0}$, where $V^{\delta _0}$ is a $\delta _0$
enlargement of $V$, \par $V^{\delta _0} = \{ z\in {\cal A}_r^l: dist
(z,V)\le \delta _0 \} $, $ ~ dist (z,V) = \inf_{y\in V} |z-y|$.
\par By virtue of Lemma 2 in \cite{whitntrams34} and Formulas $(1.1)$, $(1.4)$
a function $\xi $ exists such that $\xi \in C^{\infty }({\cal
A}_r^l,{\bf R})$, $ ~ \xi (z)=1$ for each $z\in V$, $\partial
^k_z\xi (z)=0$ for each $z\in V$ and $k\in {\bf N}_0^{l2^r}$, $\xi
(z)=0$ for each $z\in {\cal A}_r^l\setminus V_1$. We consider a
function $h(z)=\xi (z)g(z)$ and define also a map
$$(3.5) \quad G(z)=T \kappa ^{l2^r} \int_{{\cal A}_r^l} h(y)
\exp ( \kappa ^2 \sum_{p=1}^l  [(\mbox{}_pz - \mbox{}_py)
\sum_{k=0}^{2^r-1}i_k((\mbox{}_pz- \mbox{}_py)i_k)])\omega (dy),$$
where $\kappa >0$ is a positive constants, $z\in {\cal A}_r^l$. \par
The function $f(z-y)= ev_y(\mu )(z)$ is ${\cal A}_r$-differentiable
in $z$. This follows from the chain rule in the sense of phrases
$\mu $ and $\nu $ and $\lambda $ of the word algebra (see 2.2.1 in
\cite{ludfov} or 2.3.2 in \cite{ludancdnb}) and Formula $(3.2)$ with
$\mu $ being a composition of suitable $\nu $ and $\lambda $. In the
considered particular case $$\nu = \sum_{n=0}^{\infty
}(n!)^{-1}X^n\quad \mbox{ and }  ~ \lambda = \sum_{p=1}^l b_p [Y_p
\sum_{k=0}^{2^r-1}i_k(Y_pi_k)]$$ over the alphabet ${\cal X}  = \{
{\cal A}_r, X, Y_1,...,Y_l, E_0,...,E_l \} $ (see also Subsections
2.1-2.5 in \cite{ludkaaca2014}). The integral on the right side of
Formula $(3.5)$ and the following integral
$$\Psi (z).v = T\kappa ^{l2^r} \int_{{\cal A}_r^l} h(y) (D_z\exp ( \kappa ^2 \sum_{p=1}^l
[(\mbox{}_pz - \mbox{}_py) \sum_{k=0}^{2^r-1}i_k((\mbox{}_pz-
\mbox{}_py)i_k)])).v \omega (dy)$$ converge uniformly in $z$ on each
compact subset $W$ in ${\cal A}_r^l$ and on $W\times \{ v\in {\cal
A}_r^l: |v|\le 1 \} $ correspondingly. From Formulas $(3.3)$ and
$(3.5)$ it follows that $D_zG(z).v= \Psi (z).v$ for each $z$ and $v$
in ${\cal A}_r^l$. By virtue of Theorems 2.11, 2.15 and 3.10 in
\cite{ludfov} the function $G(z)$ is ${\cal A}_r$ analytic on ${\cal
A}_r^l$.
\par In view of the theorem about differentiation of an improper
integral depending on parameters (see Theorem 520.3 in
\cite{fihteng} or Subsection XVII.2.3 in \cite{zorich})
$$(3.6)\quad \frac{\partial ^{|k|} G(z)} {\partial
\mbox{}_1z_0^{k(1,0)} ...\partial \mbox{}_lz_{2^r-1}^{k(l,2^r-1)}} =
(-1)^{|k|} T\kappa ^{l2^r}\cdot $$ $$ \int_{{\cal A}_r^l} h(y)
\frac{\partial ^{|k|}\exp ( \kappa ^2 \sum_{p=1}^l  [(\mbox{}_pz -
\mbox{}_py) \sum_{k=0}^{2^r-1}i_k((\mbox{}_pz- \mbox{}_py)i_k)])}
{\partial \mbox{}_1y_0^{k(1,0)} ... \partial
\mbox{}_ly_{2^r-1}^{k(l,2^r-1)}}\omega (dy),$$ where
$k=(k(1,0),...,k(1,2^r-1),...,k(l,0),...,k(l,2^r-1))$;
\\ $|k|=k(1,0)+...+k(l,2^r-1)$; $ ~ k(p,j)\in {\bf N}_0$ for each
$p$ and $j$. The integration by parts transforms the right part of
Formula $(3.6)$ into
$$(3.7)\quad \frac{\partial ^{|k|} G(z)} {\partial
\mbox{}_1z_0^{k(1,0)} ...\partial \mbox{}_lz_{2^r-1}^{k(l,2^r-1)}} =
T\kappa ^{l2^r}\cdot $$ $$ \int_{{\cal A}_r^l} \frac{\partial
^{|k|}h(y)} {\partial \mbox{}_1y_0^{k(1,0)} ...
\partial \mbox{}_ly_{2^r-1}^{k(l,2^r-1)}} \exp ( \kappa ^2 \sum_{p=1}^l
[(\mbox{}_pz - \mbox{}_py) \sum_{k=0}^{2^r-1}i_k((\mbox{}_pz-
\mbox{}_py)i_k)]) \omega (dy).$$ From the choice of the function
$h(z)$ it follows that $$\sup_{y\in {\cal A}_r^l; ~ |k|\le m}
|\frac{\partial ^{|k|}h(y)} {\partial \mbox{}_1y_0^{k(1,0)} ...
\partial \mbox{}_ly_{2^r-1}^{k(l,2^r-1)}}|=:K_h<\infty $$
and for each $\epsilon >0$ a positive number $\delta $ exists so
that
$$\sup_{|z-y|<\delta ; ~ y\in {\cal A}_r^l; ~ z \in {\cal A}_r^l; ~ |k|\le m}
|\frac{\partial ^{|k|}h(z)} {\partial \mbox{}_1z_0^{k(1,0)} ...
\partial \mbox{}_lz_{2^r-1}^{k(l,2^r-1)}}- \frac{\partial ^{|k|}h(y)} {\partial \mbox{}_1y_0^{k(1,0)} ...
\partial \mbox{}_ly_{2^r-1}^{k(l,2^r-1)}}|\le \epsilon /2.$$
A positive number $\kappa _0>0$ exists such that $1-\Phi (\kappa
\delta )<\epsilon /(4K_h)$ for each $\kappa \ge \kappa _0$.
Therefore for any $z\in {\cal A}_r^l$ we get that
$$(3.8)\quad T\kappa ^{l2^r}|\int_{B({\cal A}_r^l,z,\delta )} [ \frac{\partial
^{|k|}h(y)} {\partial \mbox{}_1y_0^{k(1,0)} ...
\partial \mbox{}_ly_{2^r-1}^{k(l,2^r-1)}} - \frac{\partial
^{|k|}h(z)} {\partial \mbox{}_1z_0^{k(1,0)} ...
\partial \mbox{}_lz_{2^r-1}^{k(l,2^r-1)}} ]\cdot $$ $$ \exp ( \kappa ^2 \sum_{p=1}^l
[(\mbox{}_pz - \mbox{}_py) \sum_{k=0}^{2^r-1}i_k((\mbox{}_pz-
\mbox{}_py)i_k)]) \omega (dy)|\le 2^{-1}\epsilon \Phi (\kappa \delta
)<\epsilon /2$$ and
$$(3.9)\quad T\kappa ^{l2^r}\cdot |\int_{{\cal A}_r^l\setminus B({\cal A}_r^l,z,\delta )} [ \frac{\partial
^{|k|}h(y)} {\partial \mbox{}_1y_0^{k(1,0)} ...
\partial \mbox{}_ly_{2^r-1}^{k(l,2^r-1)}} - \frac{\partial
^{|k|}h(z)} {\partial \mbox{}_1z_0^{k(1,0)} ...
\partial \mbox{}_lz_{2^r-1}^{k(l,2^r-1)}} ]\cdot $$  $$ \exp ( \kappa ^2 \sum_{p=1}^l
[(\mbox{}_pz - \mbox{}_py) \sum_{k=0}^{2^r-1}i_k((\mbox{}_pz-
\mbox{}_py)i_k)]) \omega (dy)|\le 2K_h[1-\Phi (\kappa \delta )]
<\epsilon /2.$$ Hence
$$|\frac{\partial ^{|k|}G(z)} {\partial \mbox{}_1z_0^{k(1,0)} ...
\partial \mbox{}_lz_{2^r-1}^{k(l,2^r-1)}} - \frac{\partial
^{|k|}h(z)} {\partial \mbox{}_1z_0^{k(1,0)} ...
\partial \mbox{}_lz_{2^r-1}^{k(l,2^r-1)}} |<\epsilon $$ for each $z\in
{\cal A}_r^l$ and $|k|\le m$, since $$\frac{\partial ^{|k|}G(z)}
{\partial \mbox{}_1z_0^{k(1,0)} ...
\partial \mbox{}_lz_{2^r-1}^{k(l,2^r-1)}} - \frac{\partial
^{|k|}h(z)} {\partial \mbox{}_1z_0^{k(1,0)} ...
\partial \mbox{}_lz_{2^r-1}^{k(l,2^r-1)}} =
T\kappa ^{l2^r}\cdot $$ $$ \int_{{\cal A}_r^l} [ \frac{\partial
^{|k|}h(y)} {\partial \mbox{}_1y_0^{k(1,0)} ...
\partial \mbox{}_ly_{2^r-1}^{k(l,2^r-1)}} - \frac{\partial
^{|k|}h(z)} {\partial \mbox{}_1z_0^{k(1,0)} ...
\partial \mbox{}_lz_{2^r -1}^{k(l,2^r-1)}} ]\cdot $$ $$ \exp ( \kappa ^2 \sum_{p=1}^l
[(\mbox{}_pz - \mbox{}_py) \sum_{k=0}^{2^r-1}i_k((\mbox{}_pz-
\mbox{}_py)i_k)]) \omega (dy).$$ From the choice of the function
$h(z)$ the assertion of this lemma follows.

\par {\bf 4. Lemma.} {\it Assume that $V$ is an open subset in
${\cal A}_r^l$ such that $V = \bigcup_{a=1}^{\infty } V_a$, where
$cl (V_a)\subset V_{a+1}$ and $V_a$ is bounded and open in ${\cal
A}_r^l$ for each $a\in {\bf N}$, some $V_a$ may be void. Let also
$g\in C^m(V, {\cal A}_{r,\bf F}^m)$, $~\epsilon _a\ge \epsilon
_{a+1} > 0$ for each $a$, where $0\le m \le \infty $, $~3\le
r<\infty $. Then an ${\cal A}_r$-analytic function $G$ on $V$ exists
such that for each $|k|\le \alpha _a$ and $a\in {\bf N}$
$$(4.1)\quad |\partial ^k_z[G(z)- g(z)]|<\epsilon _a \mbox{ on }
(V\setminus V_a),$$ where $\alpha _a=m$ if $m$ is finite, $\alpha
_a=a$ if $m=\infty $.}
\par {\bf Proof.} Note that in the case $V_1=\emptyset $,...,$V_s=\emptyset $
the assertion $(4.1)$ means
\par $(4.2)\quad |\partial ^k[G(z)-g(z)]|<\epsilon _s$ \\ on $V$
for each $|k|\le \alpha _s$. Next we put $\hat{W}_a=cl (V_{a-1})$,
$~W_a=(cl (V_{a+1}))\setminus V_a$, $\check{W}_a={\cal
A}_r^l\setminus V_{a+2}$ and take the closed set $\hat{W}_a\cup
W_a\cup \check{W}_a$. In view of Lemma 2 in \cite{whitntrams34} and
Formulas $(1.1)$, $(1.4)$ for each $a\in {\bf N}$ a function $f_a\in
C^{\infty }({\cal A}_r^l,{\bf R})$ exists such that $f_a(z)=1$ for
each $z\in W_a$, also $f_a(z)=0$ for each $z\in \hat{W}_a\cup
\check{W}_a$; $\partial ^k_zf_a(z)=0$ for each $z\in \hat{W}_a\cup
W_a\cup \check{W}_a$ and $0<|k|$. Therefore \par $(4.3)$ $N_a :=
\max (1, \sup_{z\in {\cal A}_r^l; ~ |k|\le \alpha _a} |\partial
^k_zf_a(z)|) <\infty $
\\ for each $a\in {\bf N}$. Then by induction in $a=1, 2,...$ analytic
functions are defined similarly to the proof of Lemma 3
$$(4.4)\quad G_a(z)=T\kappa _a^{l2^r} \int_{{\cal A}_r^l} f_a(y) [g(y)-\sum_{s=0}^{a-1}G_s(y)] \cdot $$
$$\exp ( \kappa ^2 \sum_{p=1}^l
[(\mbox{}_pz - \mbox{}_py) \sum_{k=0}^{2^r-1}i_k((\mbox{}_pz-
\mbox{}_py)i_k)]) \omega (dy) ,$$ since $f(y)\in {\bf R}$ and the
real field $\bf R$ is the center of the Cayley-Dickson algebra
${\cal A}_r$, where $3\le r<\infty $, $~G_0=0$. For each $a\in {\bf
N}$ a constant $0<\kappa _a<\infty $ exists such that
\par $(4.5)$ $|\partial_z^k(G_a(z)-H_a(z))|< \nu _a$ \\ on $cl (V_{a+1})$ for
each $k\in {\bf N}_0^n$ with $|k|\le \alpha _{a+1}$, where
$H_a(z)=f_a(z)[g(z)-G_1(z)-...-G_{a-1}(z)]$ for all $z$, $ ~ \nu _a=
\epsilon _{a+1} 2^{-a-2} [(\alpha _{a+1}+1)!]^{-n}N_{a+1}$, where
$n=l2^r$. From Lemma 3 it follows that
\par $(4.6)$ $|\partial_z^k[g(z)-G_1(z)-...-G_{a-1}(z)]|<\epsilon _a2^{-a-2}$ \\ on
$W_a$ for each $|k|\le \alpha _{a+1}$, since $\nu _a<\epsilon
_a2^{-a-2}$. From $(4.3)$ and $(4.6)$ we deduce that $|\partial
^k_zH_a(z)|<\epsilon _a2^{-a-1}$ on $W_{a-1}\cup V_{a-1}$ for each
$|k|\le \alpha _a$. From the latter inequality and $(4.5)$ it
follows that \par $(4.7)$ $|\partial _z^kG_a(z)|<\epsilon _a2^{-a}$
\\ on $cl (V_a)$ for each $|k|\le \alpha _a$, since $f_a(z)=1$ for
all $z\in V_a$. \par Put $G(z)=\sum_{j=1}^{\infty } G_j(z)$. Then
the series $\sum_{j=1}^{\infty } \partial ^k_zG_j(z)$ converges
uniformly on any compact subset of $V$, hence $\partial ^k_zG(z)=
\sum_{j=1}^{\infty } \partial ^k_zG_j(z)$ on $V$ for each $|k|\le m$
if $m$ is finite or $|k|<\infty $ if $m=\infty $. Moreover, together
with $(4.7)$ this implies that $|\sum_{j=1}^{\infty } \partial
^k_zG_{a+j}(z)|<\epsilon _a/2$ on $cl (V_{a+1})$ for each $|k|\le
\alpha _{a+1}$. Therefore, from $(4.6)$ it follows that $|\partial
_z^k (G(z)-g(z))|<\epsilon _a$ on $V_a$ for each $a\in \bf N$ and
$|k|\le \alpha _a$. Thus $(4.1)$ is proven.
\par It remains to verify that the function $G(z)$ is analytic on $V$.
Take an arbitrary point $\xi $ in $V$ and choose $a_0\in {\bf N}$
sufficiently large such that $2^{-a_0/2}<\rho <\infty $ and $B({\cal
A}_r^l,\xi ,\rho )\subset V_{a_0}$, since $cl (V_a)\subset V_{a+1}$
for each $a$ and $\bigcup_{a=1}^{\infty }V_a=V$. Put $m_a=\max_z
|H_a(z)|$. Then from $(4.4)$ we infer that
$$(4.8)\quad |(D_zG_a(z)).u|<T\kappa _a^{l2^r} m_a \int_{V_{a+2}\setminus
V_{a-1}} 2\kappa _a^2 \rho \cdot \exp (- \kappa _a^2 \rho ^2) |u|
\omega (dz)$$ $$< 2T\rho \kappa _a^{2+l2^r} m_a \exp (- \kappa _a^2
\rho ^2) v_{a+2}|u| ~ \mbox{ and}$$
$$(4.9)\quad |G_a(z))|<T\kappa _a^{l2^r} m_a \int_{V_{a+2}\setminus
V_{a-1}} \exp (- \kappa _a^2 \rho ^2) \omega (dz)$$
$$< T \kappa _a^{l2^r} m_a \exp (- \kappa _a^2 \rho ^2)
v_{a+2}$$ for each $z\in V$ and $a>a_0$ and $u\in {\cal A}_r^l$,
since $|bc|\le |b||c|$ for each $b$ and $c$ in ${\cal A}_r$, where
$v_a=\int_{V_a}\omega (dz)$. Let by induction in $a\in \bf N$ a
constant $0<\kappa _a<\infty $ be chosen such that $|G_a(z)|<2^{-a}$
on $B({\cal A}_r^l,\xi ,\rho )$ and  $|(D_zG_a(z)).u|<2^{-a}|u|$ on
$B({\cal A}_r^l,\xi ,\rho )\times {\cal A}_r^l$ and $(4.5)$ be valid
on $cl(V_{a+1})$ for each $k\in {\bf N}_0^n$ with $|k|\le \alpha
_a$, where $z\in B({\cal A}_r^l,\xi ,\rho )$, $u\in {\cal A}_r^l$.
Hence the series $\sum_{j=1}^{\infty }G_j(z)$ and
$\sum_{j=1}^{\infty } D_zG_j(z).u$ uniformly converge on $B({\cal
A}_r^l,\xi ,\rho )$ and on $B({\cal A}_r^l,\xi ,\rho )\times B({\cal
A}_r^l,0,1)$ respectively, where $z\in B({\cal A}_r^l,\xi ,\rho )$
and $u\in B({\cal A}_r^l,0,1)$. By virtue of theorems 2.11, 2.15 and
3.10 in \cite{ludfov} the function $G(z)$ is ${\cal A}_r$-analytic
on $V$.

\par {\bf 5. Lemma.} {\it Assume that $V$, $m$, $f_a$, $W_a$ and $\epsilon _a$
are provided as in  Lemma 4, also a monotonously non-decreasing
sequence $\eta _a$ of positive continuous functions $\eta _a :(0,
\infty ) \to (0, \infty )$ is defined having the limit
$\lim_{t\downarrow 0} \eta _a(t)=0$ for each $a\in \bf N$, $~w\in V$
is a marked point, $M>0$ is a positive constant. Then a sequence $
\{ \kappa _a: ~ a\in {\bf N} \} $ of positive numbers exists such
that from the conditions $(5.1)-(5.4)$:
\par $(5.1)$ $g\in C^m(V,{\cal A}_{r,\bf F})$ and \par $(5.2)$
$|g(w)|\le M$ and \par $(5.3)$ $|\partial ^k_zg(z)-\partial ^k_y
g(y)|<\eta _a(|z-y|)$ on $cl (V_a)$ for each $|k|\le \alpha _a$ and
$a\in {\bf N}$ and \par $(5.4)$ $G(z)$ is defined through $g$ as in
Lemma 4, \par it follows that $G$ is ${\cal A}_r$-analytic on $V$
and Formula $(4.1)$ is accomplished.}
\par {\bf Proof.} For each $a\in {\bf N}$ a continuous function
$\mu _a: [0,\infty )\to [0, \infty )$ exists possessing properties
$(5.5)-(5.7)$:
\par $(5.5)$ $|\partial ^k_zf_a(z) - \partial ^k_yf_a(y)|<\mu
_a(|z-y|)$ for every $y$ and $z$ in ${\cal A}_r^l$ and $|k|\le
\alpha _a$ and $a\in \bf N$,
\par $(5.6)$ $\mu _a(t)\le \mu _{a+1}(t)$ for each $a\in {\bf N}$
and $t\in [0, \infty )$,
\par $(5.7)$ $\lim_{t\downarrow 0} \mu _a(t)=0$ for each $a\in {\bf
N}$, \\ since $f_a(z)$ and $\partial ^k_zf_a(z)$ are uniformly
continuous on ${\cal A}_r^l$ for every $a\in {\bf N}$ and $k\in {\bf
N}_0^{l2^r}$. From the conditions imposed on $g$ it follows that a
constant $0<\beta _1<\infty $ exists such that for each $z\in cl
(V_3)$ and $|k|\le \alpha _2$ the inequality is valid $|\partial
^k_z g(z)|<\beta _1$. Moreover $\sup_{z\in {\cal A}_r^l: ~ |k|\le
\alpha _{a+1}} |\partial ^k_zf_a(z)|=: \gamma _a <\infty $ for each
$a\in \bf N$. On the other hand, $f_1|_{V\setminus V_3}=0$,
consequently,
\par $|\partial ^k_zf_1(z)\partial ^s_zg(z) - \partial
^k_yf_1(y)\partial
^s_yg(y)|\le \beta _1 \mu _1(|z-y|)+ \gamma _1 \eta _3(|z-y|)$ \\
for each $|k|\le \alpha _2$ and $|s|\le \alpha _2$. Therefore for
$a=1$
\par $(5.8)$ $|\partial ^k_zH_a(z) -  \partial ^k_yH_a(y)|\le \zeta
_a(|z-y|)$ \\ for every $y$ and $z$ in ${\cal A}_r^l$ and $|k|\le
\alpha _{a+1}$, where $\zeta _a(t) := [(\alpha
_{a+1}+1)!]^{l2^r}[\beta _a \mu _a(t)+ \gamma _a \eta _{a+2}(t)$ for
each $0\le t$, where $H_1=f_1g$. Moreover, for $a=1$ we get
$|\partial ^k_zH_a(z)|\le [(\alpha _{a+1}+1)!]^{l2^r}\beta _a \gamma
_a$ for every $z\in cl (V_{a+2})$ and $|k|\le \alpha _{a+1}$, also
$\partial ^k_zH_a(z)|_{{\cal A}_r^l\setminus V_{a+2}}=0$ for each
$|k|\le \alpha _{a+1}$.
\par Hence for $a=1$ there exists $0<\delta _a<\infty $ so that $\zeta _a(t)< \nu _a2^{-a} $
for each $0\le t<\delta _a$. Then for $a=1$ we choose $0<\kappa
_{a,0}<\infty $ in such a way that $1- \Phi (\kappa _a\delta _a) <
\nu _a[4[(\alpha _{a+1}+1)!]^{l2^r}\beta _a\gamma _a]^{-1}$ for each
$\infty > \kappa _a\ge \kappa _{a,0}$ (see Formula $(3.4)$). For an
admissible function $g$ a function $G_1$ is prescribed by Formula
$(4.4)$. Hence Formula $(4.5)$ is satisfied with $a=1$. Then for
$a=1$ and an arbitrary fixed point $\xi \in V$ we take $\kappa
_{a,0}<\kappa _a<\infty $ so that for each $z\in B({\cal A}_r^l,\xi
,\rho _a)$
$$(5.9)\quad |(D_zG_a(z)).u|<
2T\rho \kappa _a^{2+l2^r} m_a \exp (- \kappa _a^2 \rho ^2)
v_{a+2}|u|<|u|2^{-a} ~ \mbox{  and}$$
$$(5.10)\quad |G_a(z)|<
T \kappa _a^{l2^r} m_a \exp (- \kappa _a^2 \rho ^2) v_{a+2}<2^{-a}$$
for each $u\in {\cal A}_r^l$, where a radius $0<\rho _a=\rho <\infty
$ is such that $B({\cal A}_r^l,\xi ,\rho _a)\subset V_a$ (see also
Formulas $(4.8)$ and $(4.9)$). Then from $(4.4)$ we deduce for $a=1$
that
$$(5.11)\quad \partial ^kG_a(z) - \partial ^k_xG_a(x)= T\kappa _a^{l2^r}
\cdot $$ $$\int_{{\cal A}_r^l} [ \partial ^k_s H_a(s)|_{s=y+z-x} -
\partial ^k_y H_a(y)] \exp ( \kappa _a^2
\sum_{p=1}^l [(\mbox{}_px - \mbox{}_py)
\sum_{k=0}^{2^r-1}i_k((\mbox{}_px- \mbox{}_py)i_k)]) \omega (dy).$$
Then from Formulas $(5.8)$, $(5.11)$, $(3.4)$ and $\lim_{t\to \infty
}\Phi (t)=1$ it follows for $a=1$ that
$$(5.12)\quad |\partial ^k_z G_a(z)- \partial ^k_x G_a(x)|\le \zeta _a(|x-z|)$$
for every $x$ and $z$ in ${\cal A}_r^l$ and $|k|\le \alpha _{a+1}$.
Then by induction in $a=1, 2, 3,...$ functions $\zeta _a$ are
defined and constants $\kappa _a$ are chosen such that the estimate
$(5.12)$ is valid and Formulas $(4.5)$, $(4.8)$, $(4.9)$, $(5.9)$
and $(5.10)$ take place for each $a<\varpi $. Then we put $H_{\varpi
}(z) = f_{\varpi }(z) [g(z) - G_1(z) -...-G_{\varpi -1}(z)]$ and
find a function $\zeta _{\varpi }$ such that Formula $(5.8)$ is
fulfilled for $a= \varpi $. Then as above we choose $\kappa _{
\varpi }$ and get $(4.5)$, $(5.9)$, $(5.10)$ and $(5.12)$ for $a=
\varpi $. As a result due to Lemma 4 there exists an ${\cal
A}_r$-analytic function $G(z)=\sum_{a=1}^{\infty }G_a(z)$ on $V$
having the property $(4.1)$.

\par {\bf 6. Theorem.} {\it Suppose that $V$ is a closed subset in ${\cal A}_r^l$
and a function $f: V\to {\cal A}_{r,\bf F}$ is of class $C^m$ in
terms of functions $f_k$, $k\in {\bf N}_0^{l2^r}$ with $|k|<m+1$,
where $0\le m\le \infty $, $l\in {\bf N}$, $3\le r<\infty $ (see
Definition 2). Then a function $F$ exists satisfying the following
conditions $(6.1)-(6.3)$: \par $(6.1)$ $F\in C^m({\cal A}_r^l, {\cal
A}_{r,\bf F})$ and \par $(6.2)$ $\partial ^k_zF(z)=f_k(z)$ for each
$z\in V$ and each either $|k|\le m$ if $m<\infty $ or $|k|<\infty $
if $m=\infty $;
\par $(6.3)$ $F(z)$ is ${\cal A}_r$-analytic on ${\cal
A}_r^l\setminus V$.}
\par {\bf Proof.} The function $f$ possesses the decomposition
\par $(6.4)$ $f=\pi _0(f)i_0+...+\pi _{2^r-1}(f)i_{2^r-1}$, \\ where $\pi _j(f):
V\to {\bf F}$, where either ${\bf F}={\bf R}$ or ${\bf F}={\bf
C}={\bf R}\oplus {\bf R}{\bf i}$. Evidently, $f$ is of class $C^m$
in terms of $f_k$ if and only if $\pi _j(f)$ is of class $C^m$ in
terms of $\pi _j(f_k)$ for each $j$. In view of Lemma 2 in
\cite{whitntrams34} applied to $A=\pi (V)$ and Formulas
$(1.1)-(1.4)$ using $z=\hat{z}(x)\in V$ with $x\in A$ an extension
$g$ of $f$ exists such that it fulfills the conditions $(6.1)$ and
$(6.2)$ (with $g$ instead of $F$) and $g\in C^{\infty }({\cal
A}_r^l, {\cal A}_{r,\bf F})$. \par We consider the complement
$\hat{V}={\cal A}_r^l\setminus V$ of $V$. Then similarly to the
proof of Lemma 4 $W_a$, $\alpha _a$ and $\epsilon _a$ are defined.
Let $F(z)=G(z)$ for each $z$ in $\hat{V}$, where $G$ is provided by
Lemma 4. Let also $F(z)=f(z)$ for each $z\in V$. Therefore, $F\in
C^m({\cal A}_r^l, {\cal A}_{r,\bf F})$ and the condition $(6.2)$ is
satisfied due to $(4.1)$. From the ${\cal A}_r$-analyticity of $G$
on $\hat{V}$ the property $(6.3)$ follows.

\par {\bf 7. Theorem.} {\it Assume that $A$ is a closed subset in ${\cal A}_r^l$
and $A_2 = \{ c_p: ~ p \in \Lambda \} $ is a subset of isolated
points in $A$, $cl (A_2)=A_2$, and $A_1= A \setminus A_2$, where
$3\le r<\infty $, $l\in {\bf N}$, $\Lambda \subset {\bf N}$.  Assume
also that either $-1\le m\in {\bf Z}$ or $m= \infty $, $ ~ m_p\in
{\bf N}_0$ for each $p\in \Lambda $ such that $$(7.1)\quad
\underline{\lim}_{p\to \infty } m_{p}\ge m$$ if $\Lambda $ is
infinite. Suppose that for each $k\in {\bf N}_0^{l2^r}$ with
$|k|<m+1$ a function $f_k(z)$ is defined on $A_1$ and for each
$|k|\le m_p$ at $c_p$, and $f$ is ${\cal A}_{r,\bf F}$ valued and of
class $C^m$ on $A$ in terms of functions $f_k$, where $f_k(c_p)=0$
for for each $m_p<|k|$ if $m_p< m$. Then a function $F\in C^m({\cal
A}_r^l, {\cal A}_{r,{\bf F}})$ exists satisfying the following
conditions $(7.2)-(7.4)$:
\par $(7.2)$ $F(z)=f(z)$ on $A$;
\par $(7.3)$ $\partial ^k_z F(z)=f_k(z)$ on $A_1$ for each $k\in {\bf N}_0^{l2^r}$
with $|k|\le m$ and at each $c_p$ for $|k|\le m_p$;
\par $(7.4)$ $F(z)$ is ${\cal A}_r$-analytic on ${\cal A}_r^l\setminus A_1$.}
\par {\bf Proof.} For a family of vectors $(j;k)$ with $k\in {\bf N}_0^{l2^r}$
and $|k|\le m_j$ for each $j\in {\bf N}$ one can choose their
enumeration $\vartheta (j;k)$ such that $\{ \vartheta (j;k): ~ j\in
{\bf N}, |k|\le m_j \} = {\bf N}$ and $\vartheta (j;k)< \vartheta
(j+1;k')$ for each $j$ and $k$ and $k'$;  $\vartheta (j;k)<
\vartheta (j;k')$ for each $j$ and $|k|<|k'|$; also $\vartheta
(j;k)\ne \vartheta (j';k')$ if $(j;k)\ne (j';k')$. Put $\theta (j;k)
= 1+\max (m_1,...,m_j,\vartheta (j;k))$ for each $j\in {\bf N}$ and
$|k|\le m_j$. For each $j\in {\bf N}$ we put \par $(7.5)$ $W_j = \{
z\in W: ~ |z|<j, d(z,A_1)>1/j \} $,
\\ where $d(z,A_1) = \inf \{ |z-y|:  ~ y \in {A_1} \} $, $W={\cal
A}_r^l\setminus A_1$. Therefore $W_j$ is open in ${\cal A}_r^l$ and
$cl (W_j)\subset W_{j+1}$ for each $j\in {\bf N}$. Then
$\bigcup_{j=1}^{\infty } W_j=W$. \par For each $j$ and $p$ in ${\bf
N}$ a radius $0<\rho (j,p)<(j+2)^{-2}$ exists such that $B({\cal
A}_r^l,c_p,\rho (j,p))\cap B({\cal A}_r^l,c_s,\rho (j,s))=\emptyset
$ for each $s\ne p$ and $B({\cal A}_r^l,c_p,\rho (j,p))\subset
W\setminus W_{\theta (p;k)}$ for each $k\in {\bf N}_0^{l2^r}$ with
$|k|\le m_p$, since $0\le m_p<\infty $ and $A_2$ is the closed set
of isolated points in $A$. By virtue of Lemma 2 in
\cite{whitntrams34} and Formulas $(1.1)-(1.3)$ and $(6.4)$ a
function $\omega _{p,k}\in C^{\infty }({\cal A}_r^l,{\cal A}_{r,{\bf
F}})$ exists for each $p\in {\bf N}$ and $|k|\le m_p$ such that
$\partial ^k_z\omega _{p,k}(z)|_{z=c_p}=1$ and $\partial ^n_z\omega
_{p,k}(z)|_{z=c_p}=0$ for each $n\ne k$ and with the restriction
$\omega _{p,k}|_{{\cal A}_r^l\setminus B({\cal A}_r^l,c_p,\rho
(j,p))}=0$, where $j$ is such that $\theta (j,k)\ge p$ and hence
$\theta (j_1,u)\ge p$ for each either $j_1>j$ or $|u|>|k|$.
Therefore a sequence $\beta _p$ exists such that $0<\beta _{p+1}\le
\beta _p<p^{-1}$ and
$$(7.6)\quad \beta _p \sup_{z\in {\cal A}_r^l} |\partial ^n_z\omega _{p,k}(z)|<
(p[(m_p+1)!]^{l2^r})^{-1}$$ for each $p\in \bf N$ and $n\in {\bf
N}_0^{l2^r}$ and $k\in {\bf N}_0^{l2^r}$ with $|n|\le m_p$ and
$|k|\le m_p$.
\par From Lemma 4 above with $\alpha _a=a$ for each $a\in \bf N$ it follows that
an ${\cal A}_r$-analytic function $H$ on $W$ exists such that
\par $|\partial ^k_z(H(z)-h(z))|<\beta _a$ for every $a\in \bf N$ and
$z\in W\setminus W_a$ and $|k|\le a$, where $h\in C^m({\cal
A}_r^l,{\cal A}_{r,\bf F})$ is an extension of $f$ with $h|_W \in
C^{\infty }(W,{\cal A}_{r,\bf F})$ and $\partial ^k_zh(z)=f_k(z)$
for every $z\in A_1$ and $k\in {\bf N}_0^{l2^r}$ with $|k|<m+1$ and
each $z=c_p$ and $|k|\le m_p$. Put $H(z)=f(z)$ for each $z\in A_1$.
Hence $H\in C^m({\cal A}_r^l,{\cal A}_{r,\bf F})$ and $\partial
^k_zH(z)=f_k(z)$ for every $z\in A_1$ and $k\in {\bf N}_0^{l2^r}$
with $|k|<m+1$ according to Theorem 6.
\par Then the choices of $\beta _p$ and $W_p$ imply that
$|b_{p,k}|<\beta _{\theta (p,k)}<\beta _p$ for each $|k|\le m_p$,
where $b_{p,k}=\partial _z^k(h(z)-H(z))|_{z=c_p}$, since $c_p\in
W\setminus W_{\theta (p,k)}$. Therefore, similarly to Lemma 9 in
\cite{whitntrams34} and Lemma 5 above a function $G: {\cal A}_r^l\to
{\cal A}_{r,\bf F}$ exists possessing the properties: its
restriction $G|_W$ is ${\cal A}_r$-analytic, $G|_{A_1}=0$, $ ~ G\in
C^m({\cal A}_r^l, {\cal A}_{r, \bf F})$ and $\partial
^kG(z)|_{z=c_p} =b_{p,k}$ for every $p\in {\bf N}$ and $|k|\le m_p$
and $\partial ^k_zG(z)=0$ for each $z\in {\cal A}_r^l\setminus W$.
Thus $F=H+G$ is the required function.

\par {\bf 8. Theorem.} {\it Suppose that $A$ and $A_s$ and $\Omega _{s+1}$
are closed subsets in ${\cal A}_r^l$ such that $A_s\subset
A_{s+1}\subset A$ and $\Omega _{s+1}\subset A_{s+1}\setminus A_s$
and $\Omega _{s+1} = \{ \hat{c}_{s+1,p}: ~ p \in \Lambda _{s+1} \}
$, $\Omega = \bigcup_s \Omega _s$ consists of isolated points
$\hat{c}_{s+1,p}$, $ ~ cl (\Omega )=\Omega $, and $\Lambda
_{s+1}\subset \bf N$ for each $-1\le s\in {\bf Z}$, where $3\le
r<\infty $ and $l\in \bf N$. Suppose also that $V\subset {\cal
A}_r^l\setminus A$ with $V=\bigcup_{s=0}^{\infty } V_s$ and
$V_s\subset  V_{s+1}$ and $cl (V_s\setminus V_{s-1})\subset (V\cup
A_s)\setminus V_{s-1}$ and $(A\cup V)\setminus V_s$ is closed for
each $-1\le s\in {\bf Z}$, where $V_{-1}=\emptyset $. Assume that a
function $f=f_0$ is ${\cal A}_{r,\bf F}$ valued and of class $C^s$
on $\Upsilon _s = (A\cup V)\setminus (A_{s-1}\cup V_{s-1})$ in terms
of functions $f_k$ with $|k|\le s$ for each $s\in {\bf N}_0$, also
that $\sigma (z)$ is a continuous function on ${\cal A}_r^l$ with
$\sigma (z)>0$ for each $z\in {\cal A}_r^l\setminus J$ and $\sigma
(z)=0$ for each $z\in J$, where $J=A\setminus \Omega $. Then a
function $F: {\cal A}_r^l\setminus A_{-1} \to {\cal A}_{r,\bf F}$
exists satisfying the following conditions: \par $(8.1)$ $F$ is of
class $C^s$ on ${\cal A}_r^l\setminus A_{s-1}$ for each $s\in {\bf
N}_0$;
\par $(8.2)$ $\partial ^k_zF(z) = f_k(z)$ for every $z\in A\setminus A_{s-1}$
and $|k|\le s$ and $s\in {\bf N}_0$, \par $(8.3)$ $|\partial
^k_zF(z)-f_k(z)|<\sigma (z)$ on $V\setminus V_{s-1}$ for each
$|k|\le s$ and $s\in {\bf N}_0$, \par $(8.4)$ $F$ is ${\cal
A}_r$-analytic on ${\cal A}_r^l\setminus J$.}
\par {\bf Proof.} At first relative to its norm ${\cal
A}_r^l$ is divided into canonically closed cubes $Q_{p,0,a}$ having
a length of ribs $t_0=1$. Consider the set $K_{p,0}$ of all such
cubes with $d(Q_{p,0,a}, A\cup R_p)\ge 6l^{1/2}2^{r-1}$, where
$d(Q,S) := \inf_{y\in Q, z \in S} |y-z|$ for subsets $Q$ and $S$ in
${\cal A}_r^l$, where $R_0=V$, $~R_p=\Upsilon _p$ for each $p\in
{\bf N}= \{ 1,2,3,... \}$. Then by induction in $s$ if sets
$K_{p,0}$, ...,$K_{p,s-1}$ of cubes are constructed each cube $Q$
which is not in $K_{p,0}\cup ... \cup K_{p,s-1}$ is divided into
cubes $Q_{p,s,a}$ having a length of ribs $2^{-s}$. Then $K_{p,s}$
is defined as the set of all cubes $Q_{p,s,a}$ satisfying the
following condition $d(Q_{p,s,a}, A\cup R_p )\ge
6l^{1/2}2^{r-1}2^{-s}$. Denote by $v_{p,\iota }$ the vertices of
cubes belonging to $K_p=\bigcup_{s=0}^{\infty }K_{p,s}$ and arranged
into a sequence with $\iota \in {\bf N}$. Choose points $x_{p,j}$ in
$\Upsilon _p$ such that $|x_{p,\iota }-v_{p,\iota }| \le 2
d(v_{p,\iota },\Upsilon _p)$ for each $\iota \in {\bf N}$.
\par Next let $\dot{Q}_1 := \{ z\in {\cal A}_r^l: ~\forall n=1,...,l , ~ \forall j=0,...,2^r-1 ~ |\mbox{}_nz_j|<1
\} $ be the interior of the unit cube, where $\mbox{}_nz\in {\cal
A}_r$, $~ \mbox{}_nz=\mbox{}_nz_0i_0+...+\mbox{}_nz_{2^r-1}
i_{2^r-1}$, $ ~\mbox{}_nz_j\in {\bf R}$ for each $n$ and $j$. We
consider the following functions at first for $p=0$:
$$(8.5)\quad \forall z\in \dot{Q}_1\setminus \{ 0 \} \quad \theta (z) = -1+ 2 \prod_{p=1}^l \prod_{j=0}^{2^r-1} (1- (\pi
_j(\mbox{}_pz))^2);$$
$$(8.6)\quad \forall z\in \dot{Q}_1\setminus \{ 0 \} \quad \theta _1(z) = \theta (z)[1-(\theta (z))^2]^{-1} ;$$
$$(8.7)\quad \forall z\in \dot{Q}_1\setminus \{ 0 \} \quad \Theta (z) = \exp (\theta _1(z));
\quad \Theta |_{{\cal A}_r^l\setminus Q_1}=0; $$
$$(8.8)\quad \forall z\in {\cal A}_r^l\setminus \{ v_{p,\iota }
\} \quad \eta _{p,\iota }(z) = \Theta (t_{p,\iota }^{-1} (
z-v_{p,\iota }));$$ $$\quad \forall z\in S_{p,\iota }\setminus
\partial S_{p,\iota } \quad \eta _{p,\iota ,1}(z)=
\Theta _1(t_{p,\iota }^{-1} ( z-v_{p,\iota }));$$
$$(8.9)\quad \forall z\in {\cal A}_r^l\setminus (A\cup V\cup Y)
\quad \phi _{p,\iota }= \eta _{p,\iota }(z)[\sum_{j} \eta
_{p,j}(z)]^{-1};$$
$$\quad \forall \iota \quad \phi _{p,\iota }(v_{p,\iota })=1; \quad \forall j\ne \iota \quad \eta _{p,\iota }(v_{p,j})=0;$$
where $t_{p,\iota }$ denotes the rib length of the largest cube in
$K_p$ having a vertex $v_{p,\iota }$, $ ~ S_{p,\iota }:= \{ z\in
{\cal A}_r^l: ~\forall n=1,...,l ~ |\mbox{}_nz-\mbox{ }_nv_{p, \iota
}|\le t_{p,\iota } \} $, $ ~ \partial S_{p,\iota }$ denotes the
boundary of $S_{p,\iota }$, $~Y= \{ v_{p,\iota }:  ~ \iota \in {\bf
N} \} $; $~ \Theta _1(z)= 1/ \Theta (z)$ for each $z\in
\dot{Q}_1\setminus \{ 0 \} $ and $\Theta _1(0)=0$. Put $$(8.10)\quad
\psi _{p;k}(z;y)= \sum_{s\in {\bf N}_0^{l2^r}; ~ |s|\le p-|k|}
\frac{f_{k+s}(y)}{s!}(\pi (z-y))^s
$$ for each $p\in {\bf N}_0$ (see Definition 2). Then a continuous
extension $g_{p,1}: {\cal A}_r^l\setminus A_{-1}\to {\cal A}_{r, \bf
F}$ of $f$ for $p=0$ exists:
$$(8.11)\quad \forall z\in {\cal A}_r^l\setminus (A\cup  V)
\quad g_{p,1}(z)=\sum_{\iota } \phi _{p, \iota } (z)\psi _{p;0}(z;
v_{p,\iota }); \quad \forall z\in \Upsilon _p \quad
g_{p,1}(z)=f(z).$$ Notice that the function $\Theta $ is ${\cal
A}_r$-analytic on $\dot{Q}_1\setminus \{ 0 \}$. In view of Theorem
7, Lemma 5 and Formulas $(8.9)$, $(8.11)$ a function $g_0: {\cal
A}_r^l\setminus A_{-1}\to {\cal A}_{r, \bf F}$ exists such that it
is ${\cal A}_r$-analytic on ${\cal A}_r^l\setminus (A\cup \Upsilon
_1)$ and
\par $g_0|_{(A\cup \Upsilon _1)\setminus A_{-1}}=g_{0,1}|_{(A\cup
\Upsilon _1)\setminus A_{-1}}$ and \par $(8.12)$ $\forall z\in {\cal
A}_r^l\setminus (A\cup \Upsilon _1)$ $ ~ |g_0(z)-g_{0,1}(z)|< \varrho _1 (z)/4$, \\
where $\varrho _s (z)=\min ( \sigma (z); d(z,(A\cup \Upsilon _s))$
for each $s\in \bf N$.
\par Then by induction in $p\in \bf N$ functions $g_p$ are
constructed such that
\par $(8.13)$ $g_p: ({\cal A}_r^l\setminus A_{-1})\to {\cal
A}_{r,\bf F}$; $ ~ g_p$ is of class $C^s$ on ${\cal A}_r^l\setminus
A_{s-1}$ for each $s=0,...,p$ and $g_p$ is ${\cal A}_r$-analytic on
${\cal A}_r^l\setminus (A\cup \Upsilon _{p+1})$;
\par $(8.14)$ $\partial ^k_zg_p(z)=f_k(z)$ for every $z\in
(A\cup \Upsilon _{p+1})\setminus A_{s-1}$ and $|k|\le s$ and
$s=0,...,p$;
\par $(8.15)$ $|\partial ^k_z(g_p(z)-g_{p-1}(z))|<\sigma
(z)2^{-p-2}$ for every $z\in (V_{p-1}\setminus V_{s-1})$ and $|k|\le
s$ and $s=0,...,p-1$;
\par $(8.16)$ $|\partial ^k_zg_p(z)- f_k(z)|<\sigma
(z)2^{-p-2}$ for every $z\in (V_p\setminus V_{p-1})$ and $|k|\le p$.
\par Let functions $g_0,...g_{p-1}$ with the properties $(8.13)-(8.16)$
be already constructed, where $p\in {\bf N}$. For ${\cal
A}_r^l\setminus (A\cup \Upsilon _p)$ a subdivision by cubes is taken
and their vertices $v_{p,\iota }$ are enumerated. Then points
$x_{p,\iota }$ are chosen in $\Upsilon _p$ so that $|x_{p,\iota
}-v_{0,\iota }| \le 2 d(v_{p,\iota },\Upsilon _p)$ (see above). Then
functions $\eta _{p,\iota }$ and $\phi _{p,\iota }$  are defined on
${\cal A}_r^l\setminus (A\cup \Upsilon _p)$ according to Formulas
$(8.8)$ and $(8.9)$. Then a function $g_{p,1}$ is defined on
$X_p\setminus \Upsilon _p$ by $(8.11)$ and put
$g_{p,1}(z)=g_{p-1}(z)$ for each $z\in \Upsilon _p$, where $X_p$ is
an open neighborhood of $\Upsilon _p$ such that $X_p\cap
A_{p-1}=\emptyset $ and $cl (X_p)\cap A_{p-1}\subset cl (\Upsilon
_p)\cap A_{p-1}$. Hence $g_{p,1}$ is a $C^p$ extension of
$f|_{\Upsilon _p}$. Therefore a function $\zeta
_p(z)=g_{p,1}(z)-g_{p-1}(z)$ for each $z\in X_p$ belongs to
$C^{p-1}(X_p,{\cal A}_{r,\bf F})$ and $\partial ^k_z\zeta _p(z)=0$
for each $z\in \Upsilon _p$ and $|k|\le p-1$.
\par Each cube in $K_p$ for each $p\in {\bf N}_0$ is compact and of finite dimension over $\bf
R$. Therefore a natural number $\chi $ exists such that for each
$p\in {\bf N}_0$ each cube $Q_{p,s,a}$ in $K_p$ has a covering
consisting at most of $\chi $ cubes $S_{p,\iota }$. Consider a set
\par $T_p := \{ z\in V_p: ~ \exists k\in {\bf N}_0^{l2^r} ~ |k|\le p
~ |\partial ^k_z\zeta
_{p+1}(z)|\ge \kappa _p \sigma (z)[\delta (z)]^{p-|k|} \} $, \\
where $\delta (z) = \max (1; d(z, \Upsilon _{p+1}))$, \par $\kappa
_p := 2^{-p-4} \chi ^{-1} [(p+1)!]^{-l2^r}(36 l^{1/2}2^{r-1})^{-p}
n_p^{-1}$ \\ for each $p\in {\bf N}_0$, $n_p=\max_{|k|\le p}N_k$,
where $N_k$ is a positive constant such that $|\partial ^k_z\phi
_{p,\iota }(z)|<2^{s|k|}N_k$ on any cube of $K_{p,s}$ for each
$\iota $. From the construction above it follows that for each $p\in
\bf N$ an open neighborhood $L_p$ of $\Upsilon _p$ exists such that
$L_p\subset X_p$ and $L_p\cap T_{p-1}=\emptyset $ and such that
$$(8.17)\quad |\partial ^k_z\zeta _p(z)|<
\kappa _{p-1}\sigma _{p,1}(z) [\delta (z)]^{p-1-|k|}\mbox{ on
}X_p,$$ where $\sigma _{p,1}(z) = \sigma (z)$ for each $z\in
V_{p-1}\cap L_p$, and $\sigma _{p,1}(z) =d(z,A_{p-1})$ for each
$z\in L_p\setminus V_{p-1}$.
\par Consider a subsequence $\{ \iota (n): ~ \forall n\in {\bf N} ~ \iota (n)<\iota (n+1) ~ S_{p, \iota (n)}
\subset L_p \} $ and put $b_p(z) = g_{p-1}(z) + \zeta _p(z) \sum_n
\phi _{p, \iota (n)}(z)$ on ${\cal A}_r^l\setminus (A\cup \Upsilon
_p)$ and $b_p(z)=f(z)$ on $(A\cup \Upsilon _p)\setminus A_{-1}$,
consequently, this function $b_p$ belongs to $C^p({\cal
A}_r^l\setminus A_{p-1}, {\cal A}_{r, \bf F})$. The Formula $(8.17)$
implies that $$(8.18)\quad |\partial ^k_z(b_p(z)-g_{p-1}(z))|<
\sigma _{p,1}(z)2^{-p-3}$$ for each $z\in ({\cal A}_r^l\setminus
A)\cup V_{p-1}$ and $|k|\le p-1$, consequently, $b_p\in C^s({\cal
A}_r^l\setminus A_{s-1},{\cal A}_{r,\bf F})$ for each $s=0,...,p$.
By virtue of Theorem 7 a function $g_p$ exists such that it is
${\cal A}_r$-analytic on ${\cal A}_r^l\setminus (A\cup \Upsilon
_{p+1})$ and $$(8.19)\quad |\partial ^k_z(g_p(z)-b_p(z))|< \varrho
_{p+1}(z)2^{-p-3}$$ for each $z\in {\cal A}_r^l\setminus (A\cup
\Upsilon _{p+1})$ and $|k|\le p$. Defining additionally
$g_p(z)=f(z)$ for each $z\in ({\cal A}_r^l\setminus A_{-1})\cup
\Upsilon _{p+1}$ provides a desired function fulfilling conditions
$(8.13)-(8.16)$. \par From the estimates $(8.18)$ and $(8.19)$ it
follows that a limit function $g(z)=\lim_p g_p(z)$ belongs to
$C^{\infty }({\cal A}_r^l\setminus A,{\cal A}_{r,\bf F})$. Then the
condition $(8.14)$ implies that $\partial ^k_zg(z)=f_k(z)$ for every
$z\in A\setminus A_{s-1}$ and $|k|\le s$ and $s\in {\bf N}_0$. Note
that $|\partial ^k_zg(z)-f_{k}(z)|< \sigma (z)2^{-1}$ for every
$z\in V\setminus V_{s-1}$ and $|k|\le s$ and $s\in {\bf N}_0$
according to $(8.11)$, $(8.12)$, $(8.15)$ and $(8.16)$.
\par The elements $\hat{c}_{s,p}$ of the two-sided ${\cal A}_r$ module
${\cal A}_r^l$ can be ordered in a sequence $c_{\iota }$. Take
$m_{\iota }=s$ if $c_{\iota }\in A_s\setminus A_{s-1}$ and consider
the set $W={\cal A}_r^l\setminus J$ and put $W_s=\{ z\in W: ~ |z|<s;
~ d(z, (A\cup V)\setminus V_{s-1})>s^{-1} \} $ for each $s\in \bf
N$. We choose $0<\rho (s,p)<(s+2)^{-2}$ a sufficiently small radius
for each $s$ and $p$ so that $B({\cal A}_r^l,c_p,\rho (s,p))\cap
B({\cal A}_r^l,c_{p_1},\rho (s,p_1))=\emptyset $ for every $p\ne
p_1$ and $s$, since $\Omega $ is the closed set of isolated points.
Therefore, for $\theta $ defined as in the proof of Theorem 7
$\forall p ~ \exists ~ k\in {\bf N}_0^{l2^r}, |k|\le m_s, ~ \exists
\theta (p,k)\ge s, ~ B({\cal A}_r^l,c_p,\rho (s,p))\subset
W\setminus W_{\theta (p,k)}$.  In view of Theorem 7 an extension $F$
of $f$ exists such that it is ${\cal A}_r$-analytic on ${\cal
A}_r^l\setminus J$, since $0<\beta _p$ can be chosen sufficiently
small by Formula $(7.6)$ so that $|\partial ^k_z(F(z)-g(z))|<\sigma
(z)/2$ for each $z\in V\setminus V_{s-1}$ and $|k|\le s$ and $s\in
{\bf N}_0$.

\par {\bf 9. Conclusion.} The results of this paper can be used
for further studies of functions of octonion and Cayley-Dickson
variables, for investigations of noncommutative geometry and
manifolds over them, for analysis and solutions of PDE which may be
with boundary conditions.


\begin{thebibliography}{199}

\bibitem{baez} J.C. Baez. "The octonions". Bull. Amer.
Mathem. Soc. {\bf 39: 2} (2002), 145-205.

\bibitem{boloktodb} N.N. Bogolubov, A.A. Logunov, A.I. Oksak,
I.T. Todorov. "General principles of quantum field theory" (Moscow:
Nauka, 1987).

\bibitem{bourbalg} N. Bourbaki. "Algebra" (Berlin: Springer, 1989).

\bibitem{dickson} L.E. Dickson. "The collected mathematical papers".
Volumes 1-5 (New York: Chelsea Publishing Co., 1975).

\bibitem{fihteng} G.M. Fihtengolz. "Course of differential and
integral calculus" V. 1-3 (Moscow: Nauka, 1966).

\bibitem{frenludkfejms18} E. Frenod, S. V. Ludkowski.
"Integral operator approach over octonions to solution of nonlinear
PDE". Far East J. of Mathem. Sci. (FJMS). {\bf 103: 5} (2018),
831-876; http://dx.doi.org/10.17654/MS103050831.

\bibitem{guetze} F. G\"ursey, C.-H. Tze. "On the role of
division, Jordan and related algebras in particle physics"
(Singapore: World Scientific Publ. Co., 1996).

\bibitem{henlei} G.M. Henkin, J. Leiterer. "Theory of
functions on complex manifolds". Series "Monographs in Mathematics";
{\bf 79} (Basel: Birkh\"auser, 1984).

\bibitem{hirschb} M.W. Hirsch. "Differential topology" (New York:
Springer-Verlag, 1976).

\bibitem{hormb3v} L. H\"ormander. The analysis of linear partial differential operators.
V. 1-4 (Springer-Verlag: Berlin, 1990).

\bibitem{hormbl} L. H\"ormander. Lectures on nonlinear hyperbolic differential equations
(Springer-Verlag: Berlin, 2003).

\bibitem{kansol} I.L. Kantor, A.S. Solodovnikov.
"Hypercomplex numbers" (Berlin: Springer-Verlag, 1989).

\bibitem{kling} W. Klingenberg. "Riemannian geometry"
(Berlin: Walter de Gruyter, 1982).

\bibitem{kodai} K. Kodaira. "Complex manifolds and
deformation of complex structures" (New York: Springer-Verlag,
1986).

\bibitem{kurosh} A.G. Kurosh. "Lectures on general algebra"
(Moscow: Nauka, 1973).

\bibitem{ludkaaca2014}  S.V. Ludkowski. "On a class of right linearly
differentiable functions of Cayley-Dickson variables". Adv. Appl.
Clifford Alg. {\bf 24: 3} (2014), 781-803.

\bibitem{ludkrimut14} S.V. Ludkowski. "Decompositions of PDE over Cayley-Dickson algebras".
Rendic. dell'Ist. di Math. dell'Univer. di Trieste. Nuova Ser. {\bf
46} (2014), 1-23.

\bibitem{ludkcvee16} S.V. Ludkowski.
"Integration of vector Sobolev type PDE over octonions". Complex
Variab. and Elliptic Equat. {\bf 61: 7} (2016), 1014-1035.

\bibitem{ludkrimut13} S.V. Ludkowski.
"Manifolds over Cayley-Dickson algebras and their immersions".
Rendic. dell'Istit. di Matem. dell'Univer. di Trieste. Nuova Ser.
{\bf 45} (2013); 11-22.

\bibitem{ludfov} S.V. Ludkovsky. "Differentiable functions of
Cayley-Dickson numbers and line integration". J. of Mathem.
Sciences, N.Y. (Springer) {\bf 141: 3} (2007), 1231-1298.

\bibitem{ludkfscdvjms07} S.V. Ludkovsky. "Functions
of several Cayley-Dickson variables and manifolds over them"; J.
Mathem. Sci.; N.Y.  {\bf 141: 3} (2007), 1299-1330.

\bibitem{ludancdnb} S.V. Ludkovsky. "Analysis over Cayley-Dickson
numbers and its applications" (Saarbr\"ucken: Lambert Academic
Publishing, 2010).

\bibitem{ludifeqcdla} S.V. Ludkovsky. "Differential equations
over octonions". Adv. Appl. Clifford Alg. {\bf 21: 4} (2011),
773-797.

\bibitem{ludkcvee13} S.V. Ludkovsky. "Integration of
vector hydrodynamical partial differential equations over
octonions".  Complex Variab. and Elliptic Equat. {\bf 58: 5} (2013);
579-609.

\bibitem{michor} P.W.  Michor.  Manifolds of differentiable
mappings (Boston: Shiva, 1980).

\bibitem{razm} Yu.P. Razmyslov. "Identities of algebras and their
representations". Series "Modern Algebra". {\bf 14} (Moscow: Nauka,
1989).

\bibitem{svalkorplb} A.G. Sveshnikov, A.B. Alshin, M.O. Korpusov,
Yu.D. Pletner. "Linear and nonlinear equations of Sobolev type"
(Moscow: Fizmatlit, 2007).

\bibitem{waerd} B.L. van der Waerden. "A history of algebra"
(Berlin: Springer-Verlag, 1985).

\bibitem{ward} J.P. Ward. "Quaternions and Cayley numbers".
Ser. Math. and its Applic. {\bf 403} (Dordrecht: Kluwer, 1997).

\bibitem{whitntrams34} H. Whitney. "Extensions of differentiable
functions". Trans. Amer. Mathem. Soc. {\bf 36} (1934), 63-89.

\bibitem{zabokoportj16} B.A. Zaikin, A.Yu. Bogadarov, A.F. Kotov, P.V. Poponov.
"Evaluation of coordinates of air target in a two-position range
measurement radar". Russian technological journal {\bf 4: 2} (2016),
65-72.

\bibitem{zorich} V.A. Zorich. "Mathematical Analysis" (Moscow: Nauka, 1984).


\end{thebibliography}
\end{document}